\documentclass[12pt,a4paper]{amsart}
\usepackage{amssymb,amsfonts,amsthm,amsmath}
\usepackage{graphicx} 

\theoremstyle{plain}
\newtheorem{theorem}{Theorem}
\newtheorem{lemma}{Lemma}

\newtheorem{corollary}{Corollary}
\newtheorem{statement}{Statement}

\numberwithin{theorem}{section}
\numberwithin{equation}{section}
\numberwithin{statement}{section}
\numberwithin{lemma}{section}
\numberwithin{definition}{section}
\numberwithin{corollary}{section} \textheight =24cm
\begin{document}
\title[vector partitions]{Polylogarithm approaches to Riemann zeta zeroes}
\author{Geoffrey B Campbell}
\address{Mathematical Sciences Institute \\
         The Australian National University \\
         ACT, 0200, Australia}

\email{Geoffrey.Campbell@anu.edu.au} 

\begin{abstract}
We use visible point vector identities to examine polylogarithms in the neighbourhood of the Riemann zeta function zeroes. New formulas limiting to the trivial zeroes and to the critical line on the zeta function are given. Similar results arise from Euler-Zagier sums given in Bailey, Borwein and Crandall ~\cite{dB2007} providing new infinite product identities.
\end{abstract}
\maketitle

\section{Polylogarithms near trivial zeroes of the Riemann zeta function} \label{S:1}

In this note we examine the limit as we approach the zeroes of the Riemann zeta function from the polylogarithm function. Although we will not hope to prove anything as profound as the Riemann Hypothesis itself, we can examine the behaviour of the polylogarithm in the vicinity of the trivial and non-trivial zeroes of the Riemann zeta function, and along the way prove some new infinite product identities. The two examples of Visible Point Vector identities cited in WolframMathWorld (see Weisstein ~\cite{eW2012}) are cases of the author's identities, which include the

\begin{theorem} \emph{(Campbell ~\cite{gC1994})} \label{T:1.1}
If $|x|, |y|, |z|<1$ and $s, t, u$, are complex numbers and $Li_s(z)= \sum_{k=1}^{\infty}\frac{z^k}{k^s}$,

\begin{equation} \label{E:1.1}
\prod_{\substack{ a,b \geq 1 \\ \gcd(a,b)=1}}
\left( 1- y^a z^b \right)^{-a^{-s}b^{-t}}
 = \exp\left[ Li_s(y) Li_t(z)\right],
  \end{equation}
\begin{center}
for which $s+t=1$.
\end{center}

\begin{equation} \label{E:1.2}
\prod_{\substack{ a,b,c \geq 1 \\ \gcd(a,b,c)=1}}
\left( 1-x^a y^b z^c \right)^{-a^{-s}b^{-t}c^{-u}}
 = \exp\left[ Li_s(x) Li_t(y) Li_u(z)\right],
  \end{equation}
\begin{center}
for which $s+t+u=1$.
\end{center}

\end{theorem}

Identities (\ref{E:1.1}) and (\ref{E:1.2}) are easy to prove by differentiating the logarithm of both sides, and applying the well known

\begin{lemma} \label{L:1.1}
Consider an infinite region raying out of the origin in any Euclidean vector space. The set of all lattice point vectors apart from the origin in that region is precisely the set of positive integer multiples of the visible point vector (vpv's) in that region. \emph{(Campbell ~\cite{gC1994})}
\end{lemma}

We can make use of some well known summations and substitute them into (\ref{E:1.1}) and (\ref{E:1.2}) above. The following cases of the polylogarithm are well known. (see Abramowitz and Stegun \cite{mA1972})

\begin{statement} \label{St:1.1}
\begin{equation} \label{E:1.3}
Li_1(z) = -\ln(1-z),
\end{equation}
\begin{equation} \label{E:1.4}
Li_0(z) = \frac{z}{1-z},
\end{equation}
\begin{equation} \label{E:1.5}
Li_{-1}(z) = \frac{z}{(1-z)^2},
\end{equation}
\begin{equation} \label{E:1.6}
Li_{-2}(z) = \frac{z(1+z)}{(1-z)^3},
\end{equation}
\begin{equation} \label{E:1.7}
Li_{-3}(z) = \frac{z(1+4z+z^2)}{(1-z)^4},
\end{equation}
\begin{equation} \label{E:1.8}
Li_{-4}(z) = \frac{z(1+z)(1+10z+z^2)}{(1-z)^5}.
\end{equation}
\end{statement}

Some particular expressions for half-integer values of the argument $z$ are:

\begin{statement} \label{St:1.2}
If $\zeta$ is the Riemann zeta function then
\begin{equation} \label{E:1.9}
Li_{1}(1/2) = \ln 2,
\end{equation}
\begin{equation}  \label{E:1.10}
Li_{2}(1/2) = \frac{\pi^2}{12} - \frac{1}{2} (\ln 2)^2,
\end{equation}
\begin{equation} \label{E:1.11}
Li_{3}(1/2) = \frac{1}{6} (\ln2)^3 - \frac{\pi^2}{12} (\ln 2)^2 - \frac{7}{8} \zeta(3).
\end{equation}
\end{statement}

No formulas of this type are known for higher integer orders (see Lewin~\cite[page 2]{lL1981}). However, one has for instance, the result from Borwein, Borwein and Girgensohn ~\cite{dB1995} that

\begin{statement} \label{St:1.3}
If the alternating double sum $\zeta(\overline{3},\overline{1}) = \sum_{m>n>0}(-1)^{m+n}m^{-3}n^{-1}$, then
\begin{equation} \label{E:1.12}
Li_4(1/2) = \frac{\pi^4}{360} - \frac{1}{24} (\ln 2)^4 - \frac{\pi^2}{24} (\ln 2)^4 - \frac{1}{2} \zeta(\overline{3},\overline{1}).
\end{equation}

\end{statement}

We can substitute from the statements ~\ref{St:1.2} and ~\ref{St:1.3} into ~\ref{E:1.1} quite easily, to obtain the expressions for $\exp\left[Li_s(x)Li_{1-s}(y)\right]$. For the first few cases, namely $s = 1, 2, 3, 4, 5$, we have,

\begin{theorem} \label{T:1.2}
If $|x|,|y|<1$,

\begin{equation} \label{E:1.13}
\prod_{\substack{ a,b \geq 1 \\ \gcd(a,b)=1}}
\left( 1- x^a y^b \right)^{-a^{-1}}
 = \left(\frac{1}{1-x}\right)^{ \frac{y}{1-y}},
  \end{equation}

\begin{equation} \label{E:1.14}
\prod_{\substack{ a,b \geq 1 \\ \gcd(a,b)=1}}
\left( 1- x^a y^b \right)^{-a^{-2}b^1}
 = \exp \left(Li_2(x) \frac{y}{(1-y)^2}\right),
  \end{equation}

\begin{equation} \label{E:1.15}
\prod_{\substack{ a,b \geq 1 \\ \gcd(a,b)=1}}
\left( 1- x^a y^b \right)^{-a^{-3}b^2}
 = \exp \left(Li_3(x) \frac{y(1+y)}{(1-y)^3}\right),
  \end{equation}

\begin{equation} \label{E:1.16}
\prod_{\substack{ a,b \geq 1 \\ \gcd(a,b)=1}}
\left( 1- x^a y^b \right)^{-a^{-4}b^3}
 = \exp \left(Li_4(x) \frac{y(1+4y+y^2)}{(1-y)^4}\right),
  \end{equation}

\begin{equation} \label{E:1.17}
\prod_{\substack{ a,b \geq 1 \\ \gcd(a,b)=1}}
\left( 1- x^a y^b \right)^{-a^{-5}b^4}
 = \exp \left(Li_5(x) \frac{y(1+y)(1+10y+y^2)}{(1-y)^5}\right),
  \end{equation}

\end{theorem}

These results explicitly stated are all new, and will allow us to examine the behavior of the Riemann zeta function in the neighborhoods of the trivial zeroes occurring at $-2,-4,-6,...$ by observing the right sides of identities (\ref{E:1.15}), (\ref{E:1.17}), etc as $y \rightarrow 1$.

As an example limiting $Li_{-2}(y)$ as $y\rightarrow1$ to the trivial zero of $\zeta(-2)=0$ we see that (\ref{E:1.15}) is equivalent to

\begin{equation} \notag
\prod_{\substack{ a,b \geq 1 \\ \gcd(a,b)=1}}
\left( 1- x^a y^b \right)^{-a^{-3}b^2}
 = \exp \left(Li_3(x) \left(y+2^2y^2+3^2y^3+4^2y^4+5^2y^5+\cdots\right)\right),
  \end{equation}

and we can examine this for the value $y = 1-\delta$, for small positive real $\delta$ which means then we have, by (\ref{E:1.15}), that

\begin{equation} \notag
\prod_{\substack{ a,b \geq 1 \\ \gcd(a,b)=1}}
\left( 1- x^a (1-\delta)^b \right)^{-a^{-3}b^2}
 = \exp \left(Li_3(x) \frac{(1-\delta)(2-\delta)}{(\delta)^3}\right)
  \end{equation}
where both sides approach 1 as $\delta \rightarrow 0$.

However, it is easily seen that the approach here is not very illuminating, as the limit of the polylogarithm approaches infinity.  A better approach to take if we are to learn something new about the zeroes of the zeta function, may be to limit along the $s+t=1$ restriction on (\ref{E:1.1}). The following particular cases are new, and incidental to theorem \ref{T:1.2} above.

\begin{corollary} \label{C:1.1}
If $|y|<1$,
\begin{equation} \label{E:1.18}
\prod_{\substack{ a,b \geq 1 \\ \gcd(a,b)=1}}
\left( 1- y^b \right)^{-a^{-2}b^1}
 = \exp \left(\frac{\pi^2}{6} \frac{y}{(1-y)^2}\right),
  \end{equation}
\begin{equation} \label{E:1.19}
\prod_{\substack{ a,b \geq 1 \\ \gcd(a,b)=1}}
\left( 1- y^b \right)^{-a^{-3}b^2}
 = \exp \left(\zeta(3) \frac{y(1+y)}{(1-y)^3}\right),
  \end{equation}
\begin{equation} \label{E:1.20}
\prod_{\substack{ a,b \geq 1 \\ \gcd(a,b)=1}}
\left( 1- y^b \right)^{-a^{-4}b^3}
 = \exp \left(\frac{\pi^4}{90} \frac{y(1+4y+y^2)}{(1-y)^4}\right),
  \end{equation}
\begin{equation} \label{E:1.21}
\prod_{\substack{ a,b \geq 1 \\ \gcd(a,b)=1}}
\left( 1- y^b \right)^{-a^{-5}b^4}
 = \exp \left(\zeta(5) \frac{y(1+y)(1+10y+y^2)}{(1-y)^5}\right),
  \end{equation}
\end{corollary}

In a similar vein, the results from statements ~\ref{St:1.2} and ~\ref{St:1.3} yield incidentally, the

\begin{corollary} \label{C:1.2}
If $|y|<1$,
\begin{equation} \label{E:1.22}
\prod_{\substack{ a,b \geq 1 \\ \gcd(a,b)=1}}
\left( 1- \frac{y^b}{2^a} \right)^{-\frac{1}{a}}
 = 2^{\frac{y}{1-y}},
  \end{equation}
\begin{equation} \label{E:1.23}
\prod_{\substack{ a,b \geq 1 \\ \gcd(a,b)=1}}
\left( 1- \frac{y^b}{2^a} \right)^{-\frac{b}{a^2}}
 = \exp \left[\left(\frac{\pi^2}{12} - \frac{1}{2} (\ln 2)^2\right) \frac{y}{(1-y)^2}\right],
  \end{equation}
\begin{equation} \label{E:1.24}
\prod_{\substack{ a,b \geq 1 \\ \gcd(a,b)=1}}
\left( 1- \frac{y^b}{2^a} \right)^{-\frac{b^2}{a^3}}
 = \exp \left[\left(\frac{1}{6} (\ln2)^3 - \frac{\pi^2}{12} (\ln 2)^2 - \frac{7}{8} \zeta(3)\right) \frac{y(1+y)}{(1-y)^3}\right],
  \end{equation}
\begin{equation} \label{E:1.25}
\prod_{\substack{ a,b \geq 1 \\ \gcd(a,b)=1}}
\left( 1- \frac{y^b}{2^a} \right)^{-\frac{b^3}{a^4}}
= \exp \left[\left(\frac{\pi^4}{360} - \frac{1}{24} (\ln 2)^4 -
 \frac{\pi^2}{24} (\ln 2)^4 - \frac{1}{2} \zeta(\overline{3},\overline{1})\right) \frac{y(1+4y+y^2)}{(1-y)^4}\right],
  \end{equation}
\begin{center}
with $\zeta(\overline{3},\overline{1})$ as in $(\ref{E:1.12})$.
\end{center}

\end{corollary}

We remark at this point that many of the Euler-Zagier sums from Borwein et al (\cite{dB1995}) and Bailey et al (\cite{dB2007}) will provide us with more new incidental corollary infinite products. In the next section we write the corresponding identities for the more interesting neighborhoods for the nontrivial zeroes of the Riemann zeta function where $y = \frac{1}{2} + iT$ for real numbers $T$.

\section{Identities near non-trivial zeroes of the Riemann zeta function} \label{S:2}

We end by taking the obvious case of (\ref{E:1.1}) to enable observation of the polylogarithm function as the argument approaches the non-trivial zeroes of the Riemann zeta function. In a similar manner to section \ref{S:1}, it is clear that (\ref{E:1.1}) from theorem \ref{T:1.1} can
be substituted with $s = \frac{1}{2} + iT$ and $t = \frac{1}{2} - iT$, whence $s+t=1$ as required of (\ref{E:1.1}) to give:

\begin{theorem} \label{T:2.1}
If $|x|, |y|<1$; and $T$ is a real number,

\begin{equation} \label{E:2.1}
\prod_{\substack{ a,b \geq 1 \\ \gcd(a,b)=1}}
\left( \frac{1}{1- x^a y^b} \right)^{\sqrt{\frac{1}{ab}}\exp\left( iT \ln\left(\frac{b}{a}\right)\right)}
 = \exp\left[ Li_{\left(\frac{1}{2} + iT\right)}(x) Li_{\left(\frac{1}{2} - iT\right)}(y)\right].
  \end{equation}

\end{theorem}


\begin{thebibliography}{99}
\bibitem{mA1972}
ABRAMOWITZ, M., and STEGUN, I.  Handbook of Mathematical Functions, Dover Publications Inc., New York, 1972.
\bibitem{tA1976}
APOSTOL, T.  Introduction to Analytic Number Theory, Springer-Verlag, New York, 1976.
\bibitem{mB2002a}
BAAKE, M., and GRIMM, U., Combinatorial problems of (quasi-)crystallography, preprint Institut fur mathematik,
Universitat Greifswald, December 2002.
\bibitem{dB2007}
BAILEY, D. H., BORWEIN, J. M., CRANDALL, R. E. Computation and theory of extended Mordell Tornheim Witten sums, 8th Asia-pacific Complex Systems Conference, International Workshop on Operator Theory and Its Applications, www.amsi.org.au/events/archived-events/197-8th-asia-pacific-complex-systemsconference, 2007.
\bibitem{wB1935}
BAILEY, W. N. Generalized Hypergeometric Series, Cambridge Univ. Press, Cambridge, (reprinted by Stechert-Hafner, New York) 1935.
\bibitem{dB1995}
BORWEIN, D; BORWEIN J. M. and GIRGENSOHN, R. (1995) "Explicit evaluation of Euler sums" Proceedings of the Edinburgh Mathematical Society (Series 2) 38 (2): 277-294. DOI:10.1017/S0013091500019088
\bibitem{gC1992}
CAMPBELL, G. B.  Multiplicative functions over Riemann zeta function products, J. Ramanujan Soc. 7 No. 1, 1992, 52-63.
\bibitem{gC1993}
CAMPBELL, G. B.  Dirichlet summations and products over primes, Internat. J. Math. \& Math. Sci.,  Vol 16, No 2, 1993, 359-372.
\bibitem{gC1994}
CAMPBELL, G. B.  Infinite products over visible lattice points, Internat. J. Math. \& Math. Sci., Vol 17, No 4, 1994, 637-654.
\bibitem{gC1997}
CAMPBELL, G. B.  Combinatorial identities in number theory related to $q$-series and arithmetical functions, Doctor of Philosophy Thesis, School of Mathematical Sciences, The Australian National University, October 1997.
\bibitem{gG1990}
GASPER G., and  RAHMAN, M.  Basic Hypergeometric Series, Encyclopedia of Mathematics and its Applications, Vol 35, Cambridge University Press, (Cambridge - New York - Port Chester - Melbourne - Sydney), 1990.
\bibitem{mG1980}
GLASSER, M. L., and ZUCKER, I. J.  Lattice Sums,  Th. Chem.: Adv. Persp.,  Vol 5,  1980, 67-139.
\bibitem{gH1971}
HARDY, G. H., and WRIGHT, E. M.  An Introduction to the Theory of Numbers,  Oxford University Press, Clarendon, London, 1971.
\bibitem{mH1996}
HUXLEY, M. N. Area, Lattice Points and Exponential Sums, London Mathematical Society Monographs, New Series 13, Oxford Science Publications, Clarendon Press, Oxford, 1996.
\bibitem{sL1991}
LANG, S.  Number Theory III, Encyclopædia of Mathematical Sciences, Vol 60, Springer-Verlag, (Berlin - Heidelberg - New York - London - Paris - Tokyo - Hong Kong - Barcelona), 1991.
\bibitem{lL1981}
LEWIN, L. (1981). Polylogarithms and Associated Functions. New York: North-Holland. ISBN 0-444-00550-1.
\bibitem{lL1991}
LEWIN, L. (Ed.) (1991). Structural Properties of Polylogarithms. Mathematical Surveys and Monographs. 37. Providence, RI: Amer. Math. Soc.. ISBN 0-8218-1634-9.
\bibitem{bN1992}
NINHAM, B. W., GLASSER, M. L., HUGHES, B. D., and FRANKEL, N. E. M$\ddot{o}$bius, Mellin, and Mathematical Physics, Physica A, 186, 1992, 441-481.
\bibitem{rS1989}
SIVARAMAKRISHNAN, R.  Classical Theory of Arithmetic Functions, Marcel Dekker, Inc., (New York and Basel), 1989.
\bibitem{eT1951}
TITCHMARSH, E. C.  The Theory of the Riemann Zeta Function, Oxford at the Clarendon Press, 1951.
\bibitem{eW2012}
WEISSTEIN, Eric W. "Visible Point Vector Identity." From MathWorld (A Wolfram Web Resource) http://mathworld.wolfram.com/VisiblePointVectorIdentity.html, 2012.
\end{thebibliography}
\end{document}